\documentclass[12pt]{article}
\usepackage{amsmath}
\usepackage{amsfonts}
\usepackage{amssymb}

 \newtheorem{thm}{Theorem}[section]
 \newtheorem{cor}[thm]{Corollary}
 \newtheorem{Lemma}[thm]{Lemma}
 \newtheorem{prop}[thm]{Proposition}
 \numberwithin{equation}{section}

\numberwithin{T}{section}

\usepackage{amsmath}
\usepackage{amsfonts}
\usepackage{amssymb}

\newcommand{\n}{\lhd}
\newcommand{\N}{{\Bbb N}}
\newcommand{\la}{{\langle}}
\newcommand{\ra}{{\rangle}}
\newcommand{\g}{{\gamma}}%

\newcommand{\CF}{{\textsc{cf}}}
\newcommand{\BCF}{{\textsc{bcf}}}
\newcommand{\CN}{{\textsc{cn}}}
\newcommand{\BCN}{\textsc{bcn}}
\newcommand{\FC}{\textsc{fc}}
\newcommand{\PP}{\textsc{p}}
\newcommand{\AP}{{\textsc{ap}}}
\newcommand{\BP}{{\textsc{bp}}}
\newcommand{\CP}{{\textsc {cp}}}%
\newcommand{\TP}{{$\tilde{\textsc p}$}} 
\newcommand{\G}{{\Gamma}}

\newcommand{\pf}{{\noindent\bf Proof.\ \ }}
\newcommand{\Dr}{\mathop{\rm Dr}}
\newcommand{\QED}{\hfill $\square$\bigskip}
\newcommand{\im}{\mbox{\rm im}}
\newcommand{\Z}{{\Bbb Z}}

\newcommand{\qed}{\hfill $\square$}

\begin{document}

\title{Groups in which each subgroup is\\ commensurable with a normal subgroup }
\date{}
\maketitle
\vskip-2.1cm
{\centerline{\em\small  to the memory of Jim Wiegold}}
\medskip
\centerline{\em Carlo Casolo, Ulderico Dardano, Silvana Rinauro}

\medskip


\date{}

\begin{abstract}
A group $G$ is  a \CN-group if 
for each subgroup $H$ of $G$ there exists a normal subgroup $N$ of $G$ such that the index  $|HN:(H\cap N)|$ is finite. 
The class of \CN-groups contains properly the classes of core-finite
groups and that of groups in which each subgroup has finite index in a normal
subgroup. 

In the present paper it is shown that a  \CN-group whose periodic images are locally finite is finite-by-abelian-by-finite. Such groups are then described into some details by considering automorphisms of abelian groups. Finally, it is shown that if $G$ is a locally graded group with the property that the above index is bounded  independently of $H$, then $G$ is finite-by-abelian-by-finite.
\end{abstract}

\noindent {\bf Key words and phrases}:  locally finite, core-finite, subnormal,  inert,  {\textsc{cf}}-group.\\ 
\noindent {{\bf 2010 Mathematics Subject Classification}: Primary  20F24, Secondary 20F18, 20F50, 20E15} 
\vskip-5mm
\maketitle
\vskip-3cm

\section{Introduction and main results}

In a  celebrated paper, B.H.Neumann \cite{BHN} showed that for a group $G$ the property that 
 each subgroup $H$  has finite index in a normal subgroup of $G$ (i.e., $|H^G:H|$ is finite) is equivalent to the fact that $G$ has  finite derived subgroup ($G$ is {\it finite-by-abelian}).\marginpar{}

 \marginpar{}A class of groups with a dual property  was considered in \cite{BLNSW}. A group $G$ is said to be a \CF-group ({\em core-finite}) if each subgroup $H$  contains a normal subgroup of $G$ with finite index in $H$ (i.e.,  $|H:H_G|$ is finite).  As Tarski groups are  {\textsc{cf}}, a complete classification of \CF-groups seems to be much difficult. However, in  \cite{BLNSW} and  \cite{SW} it has been proved that {\it a \CF-group $G$ whose periodic quotients are locally finite is abelian-by-finite and, if $G$ is periodic,  there exists an integer $n$ such that  $|H:H_G|\le n$ for all $H\le G$} (say that {\it$G$ is \BCF, boundedly \CF}) and that {\it a locally graded \BCF-group is 
abelian-by-finite}. Furthermore, an easy example of a metabelian (and even hypercentral) group which is \CF\ but not \BCF\ is given.  It seems to be a still open question 
whether every locally graded \CF-group is abelian-by-finite. 
Recall that a group is said to be {\em  abelian-by-finite} if has an abelian subgroup with finite index and that a group is said to be {\em  locally finite} ({\em  locally graded}, respectively) if each non-trivial  finitely generated subgroup is finite (has a proper subgroup with finite index, respectively). 

With the aim of considering the above properties in a common framework, 
recall that two  subgroups $H$ and $K$ of a group $G$ are said to be {\em commensurable}
if $H\cap K$ has finite index in both $H$ and $K$. 
This is an equivalence relation and will be denoted by $\sim$. Clearly, if $H\sim K$, then $(H\cap L)\sim (K\cap L)$ and $HM\sim KM$ for each $L\le G$ and $M\n G$. 

Thus, in the present paper we consider the class of 
\emph{\CN-groups}, that is, { \em groups in which each subgroup is commensurable with a normal subgroup}. Into details, for a subgroup $H$ of a group $G$ define $\delta_G(H)$ to be the minimum index $|HN:(H\cap N)|$ with $N\n G$. Then $G$ is a \CN-group if  and only if $\delta_G(H)$ is finite for all $H\le G$. Clearly, subgroups and quotients of \CN-groups are also \CN-groups. 

Note that if a subgroup $H$ of a group $G$ is commensurable with a normal subgroup $N$, then  $S:=(H\cap N)_N$ has finite index in $H$. 
 Thus the class of \CN-groups is contained in the class of {\em sbyf-groups}, that is, groups in which each subgroup  $H$ contains a subnormal subgroup $S$ of $G$ such that the index $|H:S|$ is finite (i.e., $H$ is  {\em subnormal-by-finite}). It is known that {\em locally finite sbyf-groups are (locally nilpotent)-by-finite} (see \cite{H.SBYF}) and {\em nilpotent-by-Chernikov} (see \cite{Cas.SBYF}).  


The extension of a finite group by a \CN-group is easily seen to be a \CN-group, see Proposition \ref{CasoElementare} below. Moreover,  from Proposition 9   in \cite{DR1}  it follows that {\em for an abelian-by-finite group properties \CN\  and \CF\ are  equivalent}. However,  for each prime $p$ there is a nilpotent $p$-group with the property \CN\  which is  neither finite-by-abelian nor abelian-by-finite,  see Proposition \ref{esempio-p-gruppo}.

Our main result is the following. 

\medskip
\noindent{\bf Theorem A} {\em Let G be a \CN-group such that every periodic image  of $G$ is locally finite. Then $G$ is finite-by-abelian-by-finite.
}
\medskip

 Here by a {\em finite-by-abelian-by-finite} group we mean a group which has a finite-by-abelian subgroup  of finite index.  The proof of { Theorem A} will be given in Sect. \ref{proofs}. The strategy of the proof will be to reduce to the case when $G$ is nilpotent and then to apply techniques of nilpotent groups theory.  To this end, 
in Sect. \ref{automorfismi}, 
we will study the  action of a \CN-group on its abelian sections.  
\medskip

We will 
consider also {\em \BCN-groups}, that is, groups $G$ for which there is $n\in \N$ such that 
$\delta_G(H)\le n$ for all $H\le G$ and prove the following theorem. 

\medskip
\noindent{\bf Theorem B} {\em  Let $G$ be a finite-by-abelian-by-finite group. \\ 
i) $G$ is \CN\ if and only if it is finite-by-\CF.\\
ii) $G$ is \BCN\ if and only if it is finite-by-\BCF.
}
\medskip

It follows that if the group $G$ is periodic and finite-by-abelian-by-finite, then $G$ is \BCN\ if and only if it is \CN. Then 
we  consider non-periodic  finite-by-abelian-by-finite \BCF-groups in 
Proposition \ref{BCF}.

\medskip

 The more restrictive property \BCN\ remains treatable when we consider the wider class of locally graded groups.

\medskip
\noindent{\bf Theorem C} {\em A locally graded \BCN-group is finite-by-abelian-by-finite.
}


\bigskip

\centerline {\bf Preliminaries}
\medskip

Our { notation} is mostly standard.  For undefined terminology and basic facts we refer to \cite{R}. 
If $\Gamma$ is a group acting on a group $G$ and $H\le G$,  
 we denote $H_{\G}:=\cap_{\gamma\in\G}H^\gamma$ and
$H^{\G}:=\la H^\gamma\ |\ \gamma\in\G\ra$. 
We say that 
$H$ is {\em $\Gamma$-invariant} (or a {\em $\Gamma$-subgroup}) if $H^\G=H$. 



We first point out a sufficient condition for a group to be \CN\ (or even \BCN) and give examples of non trivial \CN-groups.

\begin{prop}\label{CasoElementare}
Let $G$ be a group with  a normal series $G_{0}\le G_{1}\le G$, where $G_{0}$ and $G/G_{1}$ have finite order, $m$ and $n$ respectively.

 If $H\le G$, then $H$ is commensurable with  $H_1:=(H\cap G_1)G_0\le G_1$ and $\delta_G(H)\le mn\cdot\delta_{G/G_0}(H_1/G_0)$.  

In particular, 
if each subgroup of $G_{1}/G_{0}$ is commensurable with a normal subgroup of  $G/G_{0}$, then $G$ is a \CN-group.\qed 
\end{prop}

\begin{prop}\label{esempio-p-gruppo}
For each prime $p$ there is a nilpotent $p$-group with the property \BCN, which is neither   abelian-by-finite nor finite-by-abelian.
\end{prop}

\pf Consider a sequence $P_n$ of isomorphic groups with order $p^4$ defined by $P_n:=\la x_n,y_n\ | \ x_n^{p^3}=y_n^p=1,\ x_n^{y_n}=x_n^{1+p^2}\ra=\la x_n\ra\rtimes \la y_n\ra$ where clearly $P_n'=\la x_n^{p^2}\ra$ has order $p$. Let $P:=\Dr_{n\in\N} P_n$ and consider the  automorphism $\gamma$ of $P$ such that $x_n^\gamma=x_n^{1+p}$ and $y_n^\gamma=y_n$, for each $n\in \N$. Clearly, $\gamma$ has order $p^2$, acts as the automorphism $x\mapsto x^{1+p}$ on $P/P'$ (which has exponent $p^2$) and acts trivially on $P'$ (which is elementary abelian). Finally let $N:=\la x_0^{p^2}x_n^{p^2}\ |\ n\in\N\ra $. Then $N$ is a $\gamma$-invariant subgroup of $P'$ with index $p$. Thus the $p$-group $G:=(P\rtimes \la\gamma\ra)/N$ is a \BCN-group by Proposition \ref{CasoElementare} applied to the series $P'/N\le P/N\le G$.

We have that  $G'$ is infinite, since for each $n$ we have $x_n^p=[x_n,\gamma]\in [P_n,\gamma]>P_n'$. Moreover, we have that $gN\in Z(P/N)$ if and only if $\forall i \ [g,P_i]\le N$, and $N\cap P_i=1$. Thus  $Z(P/N)=Z(P)/N$ where $Z(P)=\Dr_n \la x_n^p\ra$ has infinite index in $P$.

If, by contradiction,  $G$ is abelian-by-finite, then there is an abelian normal subgroup $A/N$ of $P/N$ with finite index. Then for some $m\in \N$ we have $P=AF$, where $F= \Dr_{n<m}P_n$ is a finite normal sugroup of $P$. Therefore $P/N$ is center-by-finite, a contradiction. 
\qed


\section{Automorphisms of abelian groups}\label{automorfismi}

Recall that an automorphism $\gamma$ of a group $A$ is said to be a {\it power automorphism} if $H^\g=H$ for each subgroup $H\le A$. \marginpar{} It is well-known  (see \cite{R}) that, if $A$ is an 
 abelian $p$-group, then there exists a $p$-adic integer $\alpha$ such that $a^\g=a^\alpha$ for all $a\in A$. Here $a^\alpha$ stands for $a^n$, where $n$ is any integer congruent to $\alpha$ modulo the order of $a$. 
On the other hand, a power automorphism of a non-periodic abelian group is either the identity or the inversion map.

As in \cite{DR1}, if $\Gamma$ is a group acting on an abelian group $A$,  
  we consider the following properties:
\\  
\ \PP)  \ \  $\forall H\le A\ \ H=H^\G$; \\  
\AP) \ $\forall H\le A\ \ |H:H_\G|<\infty$;\\
\BP) \ $\forall H\le A\ \ |H^\G:H|<\infty$;\\
\CP) \ $\forall H\le A\ \ \exists K=K^\G\le A$ {\em such that $H \sim K$}, \ ({\em $H$, $K$ are  commensurable}).

\smallskip
Obviously both \AP\ and \BP\ imply \CP. Moreover,  from Propositions 8 and 9 in \cite{DR1} it follows that {\em these three properties  are equivalent, 
provided $A$ is abelian and $\G$ is finitely generated, while they are in fact different in the general case even when $A$ and $\G$ are elementary abelian $p$-groups}. 
On the other hand, the properties \AP\ and \BP\ have been previously characterized in \cite{FGN} and \cite {Cas} respectively, as we are going to recall.

  To shorten statements we define a further property:
\vskip1mm 
\noindent \TP) \ 
{\em 
$\G$ has \PP\ on the factors of a $\G$-series $1\le V\le D\le A$ where\\
\phantom{m} i)\ \   $V$ is free abelian of finite rank,\\ 
\phantom{m} ii)  $D/V$ is divisible periodic with finite total rank, \\
\phantom{m}  iii) $A/D$ is periodic and has finite $p$-exponent for each prime $p\in\pi(D/V)$.  }

\begin{thm}\label{Recall1} {\rm \cite{FGN},\cite{Cas}}\  Let  $\Gamma$ be group acting on an abelian group $A$. Then:\\
a)\ $\G$ has \AP\  on $A$ if and only if there is a $\G$-subgroup $A_1$ such that
$A/A_1$ is finite and $\G$ has either \PP\ or \TP\ on $A_1$.\\ 
b)\   \ $\G$ has \BP\  on $A$ if and only if there is a $\G$-subgroup\ $A_0$\ such that
$A_0 $ is\  finite and $\G$ has either \PP\ or \TP\  on $A/A_0$.
\end{thm}

By next statement we give a characteration of the property \CP\  along the same lines.  

\begin{thm}\label{TeoremaCP}
Let  $\Gamma$ be group acting on an abelian group $A$. Then:\\ 
c) $\Gamma$ has \CP\ on $A$ if and only if  there are $\G$-subgroups $A_0\le A_1\le A$ such that $A_0$ and $A/A_1$ are finite and $\G$ has either \PP\ or \TP\  on $A_1/A_0$.
\end{thm}

\noindent The proof of Theorem \ref{TeoremaCP} is at the end of this section. Here we deduce a corollary.

\begin{cor}\label{CorollarioCP} 
For a group $\Gamma$ acting on an abelian group $A$, the following are equivalent:\\
a) $\Gamma$ has \AP\ on $A/A_0$ for a finite $\Gamma$-subgroup $A_0$ of  $A$,\\ 
b) $\Gamma$ has \BP\ on a finite index $\Gamma$-subgroup $A_1$ of  $A$, \\ 
c) $\Gamma$ has \CP\ on $A$.\qed
\end{cor}

Let us state a couple of elementary basic facts.

%

\begin{prop}\label{componenti} Let  $\Gamma$ be group acting on a locally nilpotent periodic group  $A$. Then 
 $\Gamma$ has \AP, \BP, \CP\  on $A$, respectively, if and only if $\G$ has \AP, \BP, \CP\  on finitely many primary components of $A$, respectively,  and \PP\ on all the other  ones.\phantom{...}
\end{prop}

\pf This proof uses the same argument as in Proposition 4.1 in \cite{DR2}. The sufficiency of the condition is clear once one notes that for each $H\le A$ it results $H=\Dr_p (H\cap A_p)$, where $A_p$ denotes the $p$-component of $A$. Concerning necessity, suppose  $\G$ does not have \PP\ on the primary $p$-component $A_p$ of $A$ for infinitely many primes $p$. Then for each such $p$ there is $H_p\le A_p$ which is not $\G$-invariant. We have that the subgroup generated by the $H_p$'s is not commensurable to any $\G$-subgroup.\QED

\begin{Lemma}\label{casibanali} 
Let  $\Gamma$ be a group acting on an abelian group $A$. If  $\G$ has \CP\ on $A$, then:\\
i)\ \ $\G$ has \PP\ on the largest periodic divisible subgroup of $A$. \\ 
ii) if $A$ is torsion-free, then  each $\g\in \G$ acts on  $A$ by either the identity or the inversion map.
\end{Lemma}

\pf Statement ($i$) follows from Lemma 4.3 in  \cite{DR2}. Concerning ($ii$), by Propositions 3.3 and 3.2 of \cite{DR2}  we have that there are coprime non-zero integers $n,m$ such that $a^m=(a^\gamma)^n$ for each $a\in A$. Consider $H$ such that $1 \ne H:=\la a_0 \ra\le A$. Then there is a $\G$-invariant subgroup $K$  of $A$ which is commensurable with $H$. Thus there is $r\in \N$ such that $K^r$ is a $\G$-invariant nontrivial subgroup of $H$. This forces $mn=\pm1$. \QED

%
%
%

Now we prove some lemmas. In the first one we do not require that the group $A$ is abelian.

\begin{Lemma}\label{LemmaX}
Let  $\Gamma$ be a group acting on a \FC-group $A$. If  $\G$ has \CP\ on $A$, then $\G$ has \BP\ on the subgroup  $X:=\{\,a\in A\,|\,\la a\ra^{\Gamma}\,{\rm is\ finite}\}$ of $A$.
\end{Lemma}

\pf Notice that $X$ is the set of elements $a$ of finite order of $A$ such that $|\G:C_{\G}(a)|$ is finite, so  $X$ is a locally finite $\G$-subgroup of $A$. For any $H\le X$ there is $K\le X$ such that $H\sim K=K^\G\le A$. Then there is a finite subgroup $F\le X$ such that $H\le KF$. Thus $H^\G\le KF^\G$ and $|H^\G:H|\le |F^\G|\cdot |HK:H|$ is finite.\qed

\begin{Lemma}\label{CongetturaC-}
Let  $\Gamma$ be a group acting on a $p$-group $A$ which is the direct product of cyclic groups.  If  $\G$ has \CP\ on $A$, then
 the subgroup $X:=\{\,a\in A\,|\,\la  a\ra^{\Gamma}\,{\rm is\ finite}\}$ has finite index in $A$.
\end{Lemma}

\pf Assume by contradiction that $A/X$ is infinite. 

Let us see, by elementary facts, that there is a sequence $(a_n)$ of elements of $A$ such that
\marginpar{}
\\
$1$)\phantom{x}  $\la a_n|n\in \N\ra=\Dr_{n\in\N}\la a_n\ra$,\\
$2$)\phantom{x}$A_I/ A_I \cap X$ is infinite,\hspace{-0.5mm} for each infinite subset $I$ \hspace{-0.5mm}of $\N$,\hspace{-0.5mm} where $A_I\hspace{-1mm}:=\hspace{-1mm} \la a_n|n\in I\ra$.

In fact, if $A/X$ has finite rank, it has a Pr\"ufer subgroup $Q/X$.
\marginpar{} 
 Let $Y$ be a countable subgroup of $A$  such that $Q=YX$.  By Kulikov Theorem (see \cite{R}) $Y$ is the direct product of cyclic groups, so that we may choose elements $a_n\in Y$ such that $\la a_n|n\in \N\ra=\Dr_{n\in\N}\la a_n\ra\le Y$ and $|a_nX|<|a_{n+1}X|$. The claim holds.
Similarly, if $A/X$ has infinite rank,  consider a countably infinite subgroup $Q/X$ of the socle of $A/X$. As above, 
let  $Y$ be a countable subgroup of $A$ such that $Q=YX$. Then we may choose elements $a_n\in Y$ which are independent mod $X$ and generate their direct product as claimed.

\medskip
We claim now that {\em there are sequences of infinite subsets $I_n$, $J_n$ of $\N$ and $\G$-subgroups $K_n\le A$ 
such that for each $n\in \N$:\\
$3)$ \phantom{xx} $I_n\cap J_n=\emptyset$ and $I_{n+1}\subseteq J_n$ \\
$4)$ \phantom{xx}  $K_n\sim A_{I_n}$\\ 
$5)$  \phantom{xx}  $(K_1\ldots K_i)\cap (A_{I_1}\ldots A_{I_n})\le (A_{I_1}\ldots A_{I_i})$,  $\forall i\le n$.
}
\medskip
\\
To prove the claim, proceed  by induction on $n$. 
Choose an infinite subset $I_1$ of $\N$ such that $J_1:=\N\setminus I_1$ is infinite. By \CP-property there exists $K_1=K_1^\Gamma$ commensurable with $A_{I_1}$.

Suppose we have defined $I_j$, $J_j$ and $K_j$ for $1\le j\le n$ such that 3-5 hold. Since $(K_1\ldots K_n)\sim (A_{I_1}\ldots A_{I_n})$, there is $m\in \N$ such that
\smallskip
\\  
6) \phantom{xx}  $(K_1K_2\ldots K_n)\cap A_{\N}\le  (A_{I_1}A_{I_2}\ldots A_{I_n}) \la a_1,\ldots, a_m\ra$.
\smallskip
\\  
  Let $I_{n+1}$ and $J_{n+1}$ be disjoint infinite subsets of $J_n\setminus\{1,\dots, m\}$. By \CP-property there exists $K_{n+1}=K_{n+1}^\Gamma$ commensurable with $A_{I_{n+1}}$. By the choice of $I_{n+1}$ it follows that
\smallskip
\\  
7) \phantom{xx} 
 $(K_1\ldots K_i)\cap  (A_{I_1}\ldots A_{I_{n+1}})\le (K_1\ldots K_i)\cap  (A_{I_1}\ldots A_{I_n})\ \ \forall i\le n$ 
\\  
and so (5) holds for $n+1$, as whished. The claim is now proved.

\medskip

Note that by (2) and (5) it follows that $A_{I_n}/ A_{I_n}\cap X$ is infinite for each $n\in\N$ and that also the following property holds
\smallskip
\\  
8) \phantom{xx} $(K_1K_2\ldots K_n)\cap  \bar A\le (A_{I_1}A_{I_2}\ldots A_{I_n}) \ \forall  n$, where $\bar A:=\Dr_{n\in \N}A_{I_n}$.

Now for each $n \in \N$, choose an element $b_n\in(A_{I_n}\cap K_n)\setminus X$.  Then we have  $B:=\la b_n\ |\ n\in \N\ra= \Dr_n \la b_n\ra$, where $\la b_{n}\ra^\Gamma$ is infinite and $\la b_{n}\ra^\Gamma\le K_{n}\sim  A_{I_{n}}$, so that\\
\smallskip  
9) \phantom{xx} $\la b_{n}\ra^\Gamma \cap A_{I_{n}}$ is infinite for each $n$.

Since there exists $B_0=B_0^\Gamma\sim B$, we may  take\\ 
\medskip 
-  $B_*:=(B_0\cap B)^\Gamma=
(B_*\cap B)^\Gamma
\le B^\Gamma$ where $B_*\sim B$.

Now  $B_*/(B_*\cap B)$ and $B/(B_*\cap B)$ are both finite and there is $n\in \N$ such that if $B_n:=\la b_1,\ldots,b_n\ra$  we have \\
- $(B_*\cap B)^\Gamma=B_*\le (B_*\cap B)B_n^\Gamma $\ \ and \\ 
- $B=(B_*\cap B)B_n$.

Since $b_n\in K_n$ for each $n$, we have $B_n\le \bar K_n:=K_{1}K_{2}\ldots K_{{n}}$ and  
\\ 
- $B^\Gamma= (B_*\cap B)^\Gamma B_n^\Gamma \le (B_*\cap B) B_n^\Gamma \le 
 (B_*\cap B)\bar K_n\le B\bar K_n$, so that\\ 
 - $B^\Gamma\cap \bar A\le  B\bar K_n\cap \bar A= B(\bar K_n\cap \bar A)\le  B A_{I_1}A_{I_2}\ldots A_{I_n}$ by (8) above.

Thus
\\
\medskip - \phantom{}
 $\la b_{n+1}\ra^\Gamma \cap A_{I_{n+1}}\le B^\Gamma\cap A_{I_{n+1}}\ \le  (B A_{I_1}A_{I_2}\ldots A_{I_n})\cap   A_{I_{n+1}} = \la b_{n+1}\ra$ is finite, a contradiction with (9).  \qed


\begin{Lemma}\label{CongetturaC+}
Let  $\Gamma$ be a group acting  on an abelian periodic reduced group $A$. 
If $\G$ has \CP\ on $A$, then there are $\G$-subgroups $A_0\le A_1\le A$ such that  $A_0$ and $A/A_1$ are finite and $\G$  has \PP\ on $A_1/A_0$.
\end{Lemma}

\pf By Proposition \ref{componenti} it is enough to consider the case when $A$ is a $p$-group. If  $A$ is the direct product of cyclic groups,  by Lemma \ref{CongetturaC-} we have that 
$A_1:=\{\,a\in A\,|\,\la a\ra^{\Gamma}\,{\rm is\ finite}\}$ has finite index in $A$. Further, by Lemma \ref{LemmaX},  $\G$ has \BP\ on $A_1$. Then the statement follows from  Theorem \ref{Recall1}.

Let $A$ be any reduced $p$-group and  $B_*$ be a basic subgroup of $A$. Then there is $ B=B^\Gamma\sim B_*$. Since $A/B_*$ is divisible, then the divisible radical of $A/B$ has finite index. Thus we may assume that $A/ B$ is divisible. By Kulikov Theorem (see \cite{R}), also $B$ is a direct product of cyclic groups, therefore by the above there are $\Gamma$-subgroups  $B_0\le B_1\le B$ such that  $B_0$ and $B/B_1$ are finite and $\G$  has \PP\ on $B_1/B_0$. We may assume $B_0=1$. Also, since $A/B_1$ is finite-by divisible, it is divisible-by-finite and we may assume it is divisible.

Let $\g\in \G$ and $\alpha$ be a $p$-adic integer such that $x^\g=x^\alpha$ for all $x\in B_1$. Consider the endomorphism $\g-\alpha$ of $A$ and note that $B_1\le \ker (\g-\alpha)$. Thus $A/ \ker (\g-\alpha)\simeq\im (\g-\alpha)$ is both divisible and reduced, hence trivial. It follows $\g=\alpha$ on the whole $A$.\phantom{.} \QED


\noindent{\bf Proof of Theorem \ref {TeoremaCP} }
For the sufficiency of the condition note that for any subgroup $H\le A$ we have $H\sim H\cap A_1$ and the latter is in turn commensurable with a $\G$-subgroup since $\G$ has \BP\ on  $A_1$ by Theorem \ref{Recall1}. 

Concerning necessity, we first prove the statement when $A$ is periodic.
Let $A=D\times R_1$, where $D$ is divisible and $R_1$ is reduced. Then there is a subgroup $R=R^\G\sim R_1$. Thus $DR$ and $D\cap R$ are $\G$-subgroups of $A$ with finite index and order respectively. Then we can assume $A=D\times R$. 
Let $X:=\{\,a\in A\,|\,\la a\ra^{\Gamma}\,{\rm is\ finite}\}$. Clearly $D\le X$, as $\G$ has \PP\ on $D$ by Lemma \ref{casibanali}. On the other hand, $X\cap R$ has finite index in $R$ by Lemma \ref{CongetturaC+}. It follows $A/X$ is finite and by Lemma \ref{LemmaX}  and Theorem \ref{Recall1} the statement holds.

 In the non-periodic case, note that if $V_0$ is a  free subgroup of $A$ such that $A/V_0$ is periodic, then there is $V_1=V_1^\Gamma\sim V_0$. Let $n:=|V_1/(V_0\cap V_1)|$. Thus by  applying Lemma \ref{casibanali} to the $\G$-subgroup  $V:=V_1^n$ we have \\ - {\em there is  a free abelian $\G$-subgroup $V$ such that $A/V$ is periodic and each $\g\in\Gamma$ acts on $V$ by either the identity or the inversion map}.

Suppose that $V$ has finite  rank. Consider now the action of $\G$ on the periodic group $A/V$ and apply the above.
 Then there is a series $V\le A_0\le A_1\le A$ such that  $A_0/V$ and $A/A_1$ are finite and $\G$  has either \PP\ or \TP\ on $A_1/A_0$. Since $A_0$ has finite torsion subgroup $T$ we can factor out $T$ and  assume $A_0=V$. Then $\G$  has  either \PP\ or \TP\ on $A_1$ as straightforward verification shows.

Suppose finally that  $V$ has infinite rank. Let $V_2\le V$ be such that $V/V_2$ is divisible periodic and its $p$-component has infinite rank for each prime $p$. We may assume $V:=V_2$. By the above case when $A$ is periodic,  there is a $\G$-series $V\le A_0 \le A_1\le A$ such that 
 $A_0/V$ and $A/A_1$ are finite and $\G$ has \PP\ on $A_1/A_0$. We may factor out the  torsion subgroup of $A_0$, as it is finite, and assume $A_0=V$. 
 
 Again let  $V_2\le V$ be such that $V/V_2$ is divisible periodic and its $p$-component has infinite rank for each prime $p$. 
Let $\g\in \G$ and, for each prime $p$, let $\alpha_p$ be a $p$-adic integer such that $x^\g=x^{\alpha_p}$ for all $x$ in the $p$-component  of $A_1/V$. Let $\epsilon=\pm1$ be such that $x^\g=x^{\epsilon}$ for all $x\in V$. By Lemma \ref{casibanali}, $\gamma$ has \PP\ on the maximum divisible subgroup  $D_p/V_2$ of the $p$-component of $A_1/V_2$. Thus $\alpha_p=\epsilon$ on $D_p/V_2$. Therefore 
$x^\gamma=x^\epsilon$ for all $x\in V$ and for all $x\in A_1/V$. 
We claim that $a^{\gamma}=a^\epsilon$ for each $a\in A_1$. To see this, for any $a\in A_1$ consider  $n\in\N$ such that $a^n\in V$. Then there is $v\in V$ such that  $a^{\gamma}=
a^{\epsilon}v$. Hence $a^{n\epsilon}=(a^{n})^\gamma=(a^{\gamma})^n = (a^{\epsilon}v)^n = a^{n\epsilon}v^n$. Thus $v^n=1$. Therefore, as $V$ is torsion-free, we have $v=1$, as whished.
\qed



\section{Proofs of the Theorems}\label{proofs} 

Recall that locally finite \CF-groups are  abelian-by-finite and \BCF\ (see \cite{BLNSW}).

\medskip
 \noindent{\bf Proof of Theorem B} It follows from Proposition  \ref{CasoElementare} and Proposition \ref{BCF1} below.\qed

\begin{prop} \label{BCF1} Let $G$ be an abelian-by-finite group. \\ 
i) \ if $G$ is \CN, then $G$ is \CF;\\ 
ii) if $G$ is \BCN, then $G$ is \BCF.
\end{prop}

\pf Let $A$ be a normal abelian subgroup with finite index $r$. Then each $H\le A$
has at most $r$ conjugates in $G$. If $\delta_G(H)\le n< \infty$ 
 then for each $g\in G$ we have $|H:(H\cap H^g)|\le 2\delta_G(H)\le 2n$ hence 
 $|H/H_G|\le (2n)^r$. More generally, if $H$ is any subgroup of $G$, then $|H/H_G|\le r(2n)^r$.
\QED

Let us characterize  \BCF-groups among abelian-by-finite \CF-groups.

\begin{prop} \label{BCF}  
Let $G$ be a non-periodic group with an abelian normal subgroup $A$ with finite index. Then  the following are equivalent:
\\
i) $G$ is a \BCF-group; \\
ii) $G$ is a \CF-group and there is $B\le A$ such that $B$ has finite exponent, $B\n G$ and each $g\in G$ acts by conjugation on  $A/B$ by either the identity or the inversion map.
\end{prop}

\pf Let $T$ be the torsion subgroup of $A$. By Lemma \ref{casibanali},
for each $g\in G$ there $\epsilon_g=\pm1$ such that $\g$ 
acts on $A/T$ as the automorphism $x\mapsto x^{\epsilon_g}$. Then the equivalence of (i) and (ii) holds with $B:=\langle A^{g-\epsilon_g}\ |\ g\in G\rangle$, by Theorem 3 of \cite{DR1}.\qed



To prove Theorem A, our first step is a reduction to nilpotent groups.

\begin{Lemma}\label{p.gruppo.risolubile} 
A soluble $p$-group $G$ with the property \CN\ is nilpotent-by-finite.
\end{Lemma}

\pf  By Theorem \ref{TeoremaCP}, one may refine the derived series of $G$ to a finite $G$-series $\cal S$ such that $G$ has \PP\ on each infinite factor of $\cal S$.
Recall that a $p$-group of power automorphisms of an abelian $p$-group is finite (see \cite{R}). Then the stability group $S\le G$ of the series $\cal S$, that is, the intersection of the centralizers in $G$ of the factors of the series, has finite index in $G$.  On the other hand,
by a theorem of Ph.Hall, $S$ is nilpotent. \qed
\bigskip

%

We recall now an elementary property of nilpotent groups.

\begin{Lemma}\label{G/Z(G)} Let $G$ be a nilpotent group with class $c$. If $G'$ has finite exponent $e$, then $G/Z(G)$ has finite exponent dividing $e^c$.

\end{Lemma}

\pf Argue by induction on $c$, the statement being clear for $c=1$.
Assume $c>1$ and that 
$G/Z$ has exponent dividing $ e^{c-1}$, where $Z/\gamma_c(G):=Z(G/\gamma_c(G))$. Then for all $g, x\in G$ we have $[g^{e^{c-1}},x]\in \gamma_c(G)\le G'\cap Z(G).$
Therefore  $1=[g^{e^{c-1}},x]^e=[g^{e^{c}},x],$
and $g^{e^{c}}\in Z(G)$, as claimed.\qed\medskip


Next lemma follows easily from Lemma 6 in {\rm \cite{M}}.

\begin{Lemma}\label{Moeheres}  Let $G$ be a nilpotent p-group and $N$ a normal
subgroup such that $G/N$ is an infinite elementary abelian group. If $H$ and $U$ are  finite subgroup
of $G$  such that\  ${ H}\cap U = 1$, there exists
a subgroup $V$ of $G$ such that $ U\le V$,  ${ H}\cap V = 1$ and $VN/N$ is infinite.\qed
\end{Lemma}
%
%

We deduce a technical lemma which is a tool for our pourpose.

\begin{Lemma}\label{MM}  
 Let $G$ be a nilpotent p-group and $N$ be a normal  \marginpar{} 
subgroup such that $G/N$ is an infinite elementary abelian group. If $N$ contains the \FC-center of $G$ and $G'$ is abelian with finite exponent, then there are  
subgroups $H$, $U$ of $G$ such that $H\cap U=1$, with injective maps $n\mapsto h_n\in H$ and    
$(i,n)\mapsto u_{i,n}\in [G,h_i^{-1}h_n]\cap U$, where $i,n\in \N$, $i< n$.
\end{Lemma}

%

\pf
Let
 \marginpar{} 
 us show that for each $n\in N$ there is  an $(n+1)$-uple $v_n:=(h_n, u_{0,n},u_{1,n},\dots, u_{n-1,n})$ of elements of $G$ such that:
\begin{enumerate}
\item $\{h_1,\dots, h_n\}$  is  linearly independent  modulo $N$;
\item $u_{i,n}\in[G,h_i^{-1}h_n]$\ \  $ \forall i\in\{0,\dots,n-1\}$;
\item $\{u_{j,h}\ |\ 0\le j<k\le n\} $ is  $\Z$-independent in $G'$;
\item $H_n\cap U_n=1$, where $H_n:=\la h_1,\dots,h_n\ra$ and  $U_n:=\la u_{j,h}\ |\ 0\le j<k\le n\ra$.
\end{enumerate}
Then the statement is true for $H:=\bigcup\limits_{n\in\N}H_n$ and $U:=\bigcup\limits_{n\in\N}U_n$.

Let  $h_0:=1$ and choose  $h_1\in G\setminus N$. Since $N$ contains the \FC-center $F$ of $G$, we have that $h_1$ has an infinite numbers of coniugates in $G$, hence $[G,h_1]$ is infinite and residually finite. Thus we may choose $u_{0,1}\in [G,h_1]$ such that  $\la u_{0,1}\ra \cap  \la h_1\ra=1$.

Assume then that we have defined $v_i$ for $i\le n$, that is, we have elements $h_0,\dots,h_n$, $u_{j,k}$, with $0\le j<k\le n$ such that conditions 1-4 hold.
 To define an adequate $v_{n+1}$,  
 note that  by Lemma \ref{Moeheres}  we have that 
{there exists $V_n\le G$ such that $H_n\le V_n$, $U_n\cap V_n=1$ and $V_nN/N$ is infinite}. Then choose \\ \smallskip 
 $i)  \phantom{xxx}   h_{n+1}\in V_n\setminus NU_nH_n.$\\ 
Note that $h_{n+1}\not \in FH_n\le NH_n$, so that  $\{h_1,\dots, h_{n+1}\}$ is independent mod $F$. In particular  $\forall i\in\{0,\dots,n\}$, $h_i^{-1}h_{n+1}\not\in F$, hence also $[G,h_i^{-1}h_{n+1}]$  is infinite. Since 
$G'$ is residually finite, we may recursively choose  $u_{0,n+1},\dots, u_{n,n+1}$
 such that $\forall i\in\{0,\dots,n\}$ \\ 
 $ii)$  \phantom{xxx} $\ u_{i,n}\in[G,h_i^{-1}h_n]$\\ 
  $iii)$  \phantom{xxx} $\la  u_{i,n+1} \ra \cap 
 U_n\la u_{h,n+1}\ |\ 0\le h<i\ra H_{n+1}=1$\\
Then properties 1-3 hold for $v_{n+1}$. Finally suppose there are $h\in H_n$, $u\in U_n$,  $s, t_0,\dots,t_n\in \Z$ such that\\		
  $iv)$  \phantom{xxx} $ a=hh_{n+1}^s=uu_{0,n+1}^{t_1}\cdots u_{n,n+1}^{t_n}\in H_{n+1}\cap U_{n+1}.$		\\ 
Then from $(iii)$ it follows $u_{n,n+1}^{t_n}=\ldots=u_{0,n+1}^{t_1}=1$. Hence $a=hh_{n+1}^s=u\in V_n\cap U_n=1$ and 4 holds.
\qed

\begin{Lemma}\label{G'} Let $G$ be a nilpotent $p$-group. If $G$ is \CN, then $G'$ has finite exponent.
\end{Lemma}

\pf 
If, by contradiction, $G'$ has infinite exponent, then the same happens to the abelian group  
 $G'/\gamma_3(G)$ and there is  $N$ such that $G'\ge N\ge \gamma_3(G)$ and  $G'/N$ is a Pr\"ufer group. We may assume $N=1$, that is, $G'$ itself is a Pr\"ufer group and $G'\le Z(G)$. Let us show that for any $H\le G$ we have $|H^G:H|<\infty$, hence $G'$ is finite, a contradiction. In fact we have that, by the \CN-property, there is  $K\n G$ such that $K\sim H$. Thus $H$ has finite index in $HK$ and we can also assume $H=HK$, that is, $H/H_G$ is finite. Thus, we can assume $H_G=1$ and $H\cap G'=1$, that is, $H$ is finite with order $p^n$ and $HG'$ is an abelian Chernikov group. It follows that $H$ is contained in the $n$-th socle $S$ of  $HG'\n G$, where $S$ is finite and normal in $G$, as whished.
\qed
%
%


\begin{Lemma}\label{nilpotente.periodico} Let $G$ be a nilpotent $p$-group. If $G$ is \CN, then $G$ is finite-by-abelian-by-finite.
\end{Lemma}
\pf Let $G$ be a counterexample. 
Then both $G'$ and $G/Z(G)$ are infinite. However, they have finite exponent by Lemmas \ref{G'} and \ref{G/Z(G)}.
Moreover, by Lemma \ref{LemmaX},  the  \FC-center $F$ of $G$ is finite-by-abelian. Thus $F$ has infinite index in $G$. On the other hand, $G/F$ has finite exponent, since $F\ge Z(G)$. 

Then $N:=FG^pG'$ has infinite index in $G$, otherwise the abelian group  $G/FG'$ has finite rank and finite exponent, hence it is finite. This implies that the nilpotent group $G/F$ is finite, a contradiction.

If $G'$ is abelian we are in a condition to apply Lemma \ref{MM} and get infinite  elements and subgroups $h_n\in H$, $u_{i,n}\in U$  as in that statement.  By \CN-property there is $K$ such that  $H\sim K\n G$. So that the set $\{h_n(H\cap K)\ /\ n\in\N\}$ is finite. 
Hence there is $i\in\N$ and an infinite set $I\subseteq \N \setminus\{1,\ldots,i\}$ such that for each $n\in I$ we have  $h_i^{-1}h_n\in  H\cap K$ and $u_{i,n}\in U \cap  [G, H\cap K]\le U\cap K$. Therefore $ U\cap K$  is infinite, in contradiction with $U\cap K\sim U\cap H=1$.

For the general case, proceed by induction on the nilpotency class $c>1$ of $G$ and assume that the statement is true for $G/Z(G)$ and even that this is finite-by-abelian.
Then there is a subgroup $L\le G$ such that $G/L$ is abelian and $L/Z(G)$ is finite. Thus $L'$ is finite and, by the above, $G/L'$ is  finite-by-abelian-by-finite, a contradiction.
\qed
\medskip

Let us consider now non-periodic  \CN-groups.

\begin{Lemma}\label{LemmaA}  Let $G$ be a \CN-group and $A=A(G)$ its subgroup generated by all infinite cyclic normal subgroups. Then $G/A$ is periodic, $A$ is abelian and each $g\in G$ acts on $A$ by either the identity or the inversion map, hence $|G/C_G(A)|\le 2$. 
\end{Lemma}

\pf For any $x\in G$ there is  $N\n G$ which is commensurable with $\la x\ra$. Then $n:=|N:(N\cap \la x\ra)|$ is finite. Thus $N^{n!}\le  \la x\ra$ where $N^{n!}\n G$. Hence  $G/A$ is periodic. 

It is clear that $A$ is abelian. Let $g\in G$. If $\la a\ra\n G$ and $a$ has infinite order, then there is $\epsilon_a=\pm 1$ such that $a^g=a^{\epsilon_a}$. On the other hand,
by  Lemma \ref{casibanali},  there is $\epsilon=\pm 1$ such that for each $a\in A$ there is a periodic element $t_a\in A$ such that $a^g=a^\epsilon t_a$. 
 It follows $a^{\epsilon_a-\epsilon}=t$. Therefore 
$\epsilon_a=\epsilon$ is independent of $a$, as wished.\qed
\medskip

\noindent{\bf Proof of Theorem A}. 
Recall from the Introduction that all subgroups of $G$ are subnormal-by-finite.
If $G$ is periodic, then, 
 by above quoted results in \cite{H.SBYF} and  \cite{Cas.SBYF} respectively,
 we may assume that $G$ is  locally nilpotent and soluble. 
Then, by Proposition \ref{componenti}, only finitely many primary components are non-abelian. Thus we may assume $G$ is a $p$-group and apply  Lemma \ref{p.gruppo.risolubile} and Lemma \ref{nilpotente.periodico}. It follows that $G$ is finite-by-abelian-by-finite.

To treat the general case, consider $A=A(G)$ as in Lemma \ref{LemmaA}. We may assume $A$ is central in $G$. Let $V$ be a torsion-free subgroup of $A$ such that $A/V$ is periodic. Then $G/V$ is locally finite and we may apply the above. Thus there is a series $V\le G_0\le G_1\le G$ such that $G$ acts trivially on $V$, $G_1/G_0$ is abelian, while $G_0/V$ and $G/G_1$ are finite. Then we can assume $G=G_1$ and note that the stabilizer $S$ of the series has finite index. Since $S$ is nilpotent (by Ph.Hall Theorem) we can assume that $G=S$ is nilpotent. If $T$ is the torsion subgroup of $G$, then $VT/T$ is contained in the center of $G/T$. Since all factors of the upper central series of $G/T$ are torsion-free we have $G/T$ is abelian. Thus $G'\le T\cap G_0$ is finite.
\qed
\medskip

\noindent{\bf Proof of Theorem C}. If the statement is false, by Theorem  A we may assume there is a counterexample $G$ periodic and not locally finite. Also we may assume $G$ is finitely generated and infinite. Let $R$ be the locally finite radical of $G$. By Theorem A again, $R$ is finite-by-abelian-by-finite. By Theorem B(i), there is a finite subgroup $G_0\n G$ such that $R/G_0$ is abelian-by-finite. We may assume $G_0=1$, so that $R$ is abelian-by-finite. 

We claim that $\bar G:=G/R$ has finite exponent at most $(n+1)!$ where $n$ is such that $n\ge \delta_G(H)$ for each $H\le G$ . In fact, for each $x\in \bar G$,  there is  $\bar N\n \bar G$ such that $|\bar N:(\bar N\cap \la x\ra)|\le n$. Thus $\bar N^{n!}\le  \la x\ra$ and  $\bar N^{n!}\n G$. Hence $\bar N^{n!}=1$ and $x^{n\cdot n!}=1$.

 By the positive answer (for all exponents) to the Restricted Burnside Problem,  there is a positive integer $k$ such that every finite image of $\bar G$ has order at most $k$. Since $\bar G$ is finitely generated, this means that the finite residual $\bar K$ of $\bar G$ has finite index and is finitely generated as well. Since also $\bar G$ is locally graded (see \cite{LMS}), we have $\bar K=1$ and $\bar G$ is finite. Therefore $G$ is abelian-by-finite, a contradiction.
\qed



Carlo Casolo, 
Dipartimento di Matematica “U. Dini”, Universit\`a di Firenze, Viale Morgagni 67A, I-50134 Firenze, Italy. \\ email: casolo@math.unifi.it

\bigskip
Ulderico Dardano,  Dipartimento di Matematica e
Applicazioni ``R.Caccioppoli'', Universit\`a di Napoli ``Federico
II'',  Via Cintia - Monte S. Angelo, I-80126 Napoli, Italy. \\ email: dardano@unina.it

\bigskip
Silvana Rinauro,  
Dipartimento di Matematica, Informatica ed Economia, Universit\`a della
Basilicata, Via dell'Ateneo Lucano 10 - Contrada Macchia Romana,
I-85100 Potenza, Italy.\\ email: silvana.rinauro@unibas.it

\end{document}